\documentclass[11pt,a4paper]{amsart}
\usepackage{amsmath,amsthm}
\DeclareMathOperator{\id}{id}

\DeclareMathOperator{\Diag}{Diag}
\DeclareMathOperator{\invol}{{\mathfrak{Invol}}}
\newcommand{\pointed}{\mathrm{po}}

\theoremstyle{definition}
\newtheorem{point}{}[section]
\newtheorem{remark}[point]{Remark}

\theoremstyle{plain}
\newtheorem{prop}[point]{Proposition}
\newtheorem{proposition}[point]{Proposition}
\newtheorem{lemma}[point]{Lemma}
\newtheorem{state}[point]{Statement}
\newtheorem{cor}[point]{Corollary}
\newcommand{\marginextend}[1]{ \addtolength{\oddsidemargin}{-#1}  \addtolength{\evensidemargin}{-#1}
  \addtolength{\textwidth}{#1}\addtolength{\textwidth}{#1}}
\newcommand{\updownextend}[1]{ \addtolength{\topmargin}{-#1}  \addtolength{\textheight}{#1}
\addtolength{\textheight}{#1}}
\DeclareMathOperator*{\unit}{\mathtt U}
\DeclareMathOperator{\Red}{R}
\marginextend{1.15cm}
\updownextend{0.0cm}
\title{Self-stabilization in certain infinite-dimensional matrix algebras}
\author{Gyula Lakos}
\address{Department of Geometry, E\"otv\"os University, P\'azm\'any P\'eter s.~1/C,  Budapest, H--1117, Hungary}
\email{lakos@cs.elte.hu}
\keywords{Linearization of loops, stabilization, Toeplitz algebra, Bott periodicity, locally convex algebras}
\subjclass[2000]{Primary: 19K99, Secondary:  47B35.}
\begin{document}
\begin{abstract}
Analytical tools to $K$-theory; namely, self-stabilization of rapidly decreasing matrices,
linearization of cyclic loops, and the contractibility of
the pointed stable Toeplitz algebra are discussed in terms of concrete formulas.
Adaptation to the $*$-algebra and finite perturbation categories is also considered.
The finite linearizability of algebraically finite cyclic loops is demonstrated.
\end{abstract}
\maketitle
\section{Introduction}
Learning $K$-theory, one likely encounters stabilization
of matrices, linearization of cyclic loops, and the contractibility of
the pointed stable Toeplitz algebra.
Stabilization of matrices is a fundamental feature of $K$-theory;
linearization of cyclic loops is an important method to prove complex Bott periodicity;
the Toeplitz algebra can also be used for the same purpose, but it is also a tool
to construct classifying spaces.
Although considered simple, these basic constructions are often treated in quite awkward
manners. The purpose of this paper is to show that these
topics can be discussed in a unified and simple way.
Our statements are formulated primarily in the setting of locally convex algebras.
This is not just for the sake of extreme generality but to demonstrate that concrete formulas and
maps can be very successful, without using approximations.
The main  statements of this paper are as follows:

\begin{state}[Self-stabilization] \label{stat:infstable}
Assume that
$\mathfrak A=\mathcal K_{\mathbb Z}(\mathfrak S)$,
i.~e.~the locally convex algebra of rapidly decreasing $\mathbb Z\times\mathbb Z$ matrices
over an other locally convex algebra $\mathfrak S$.
Let $r:\mathcal K_{\mathbb Z}(\mathfrak A)\rightarrow\mathfrak A$
be an isomorphism which comes from relabeling $\mathcal K_{\mathbb Z\times\mathbb Z}(\mathfrak S)$ into
$\mathcal K_{\mathbb Z}(\mathfrak S)$. Then there is a smooth homotopy
\[E:\mathcal K_{\mathbb Z}(\mathfrak A)\times[0,\pi/2]
\rightarrow\mathcal K_{\mathbb Z}(\mathfrak A)\]
such that it yields a family of endomorphisms of
$\mathcal K_{\mathbb Z}(\mathfrak A)$,
which are isomorphisms for $\theta\in[0,\pi/2)$, and a
closed injective endomorphism for $\theta=\pi/2$, with
\[E(A,0)=A,\qquad\qquad E(A,\pi/2)=\Diag(\ldots,0,0\,|\,r(A),0,0,\ldots);\]
cf. \eqref{eq:diag} for the diagonal notation.
This statement extends to unit groups, showing that $\unit(\mathcal K_{\mathbb Z}(\mathfrak A))$
can be pushed down by a homotopy into $\unit(\mathfrak A\mathbf e_{00})$.
\end{state}
(The continuous map $\phi:\mathfrak A\times[0,\pi/2]_\theta\rightarrow\mathfrak B$
is smooth in the variable $\theta$ if the higher partial derivatives
$\partial_\theta^n\phi:\mathfrak A\times[0,\pi/2]_\theta\rightarrow\mathfrak B$ are still
continuous functions.)

Let $\mathfrak A[\mathsf z^{-1},\mathsf z]$ be the algebra of formal Laurent series with rapidly decreasing
coefficients.
\begin{state}[Linearization of cyclic loops] \label{stat:linloop}
There is a smooth homotopy
\[K:\unit(\mathfrak A[\mathsf z^{-1},\mathsf z])\times[0,\pi/2]
\rightarrow\unit(\mathcal K_{\mathbb Z}(\mathfrak A)[\mathsf z^{-1},\mathsf z]),\]
such that
\[K(a(\mathsf z),\pi/2)=\Diag(\ldots1,1\,|\,a(\mathsf z)a(1)^{-1},1,1,\ldots),\]
but
\[\notag K(a(\mathsf z),0)=\mathsf U(a)\Lambda(\mathsf z,\mathsf Q)\mathsf U(a)^{-1}\Lambda(\mathsf z,\mathsf Q)^{-1};\]
where $\Lambda(\mathsf z,\mathsf Q)=\Diag(\ldots,\mathsf z,\mathsf z\,|\,1,1,\ldots)$ is the linear loop generated by the Hilbert transform, cf. \eqref{eq:diag},
and $\mathsf U(a)$ is the matrix of multiplication by $a(\mathsf z)$, cf. \eqref{eq:u}.
\end{state}

\begin{state}[Toeplitz contractibility] \label{stat:toepcont}
Let $\mathfrak A=\mathcal K_{\mathbb Z}(\mathfrak S)$.
Then the unit group of the pointed  Toeplitz algebra over  $\mathfrak A$, i. e.
 $\unit(\mathcal T_{\mathbb N}(\mathfrak A)^{\mathrm{po}})$, is contractible.
\end{state}

These statements were formulated in the smooth category.
However, it is often useful to work in slightly different categories.
One  case is when $\mathfrak A$ is a $*$-algebra. In those cases,
instead of the general unit group $\unit(\mathfrak A)$ of invertible elements,
one should work with the group $\unit^*(\mathfrak A)$ of unitary elements.
Another type of restriction occurs in the finite perturbation category, when
the algebra $\mathcal K_{\mathbb Z}(\mathfrak S)$ of rapidly decreasing matrices
 is replaced by
the algebra $\mathcal K^{\mathrm f}_{\mathbb Z}(\mathfrak S)$ of matrices with finitely many  nonzero entries,
and the algebra $\mathfrak A[\mathsf z^{-1},\mathsf z]$ of rapidly decreasing Laurent series is
replaced by the algebra $\mathfrak A[\mathsf z^{-1},\mathsf z]^{\mathrm f}$
of finite Laurent series.
(Here one should be careful, because for smooth loops being finite and invertible does
not generally imply that the inverse is finite.)

\begin{state}\label{stat:finite}
Statements \ref{stat:infstable}--\ref{stat:toepcont} restrict to the
$*$-algebra and/or finite perturbation categories.
\end{state}

The setting of finite perturbations, may, however, be too restrictive.
Let us call an element $a(\mathsf z)\in\unit(\mathfrak A[\mathsf z^{-1},\mathsf z])$ algebraically finite if $a=a_s\ldots a_1$, where for each $s$ either
$a_s$ or $(a_s)^{-1}$ has  finite Laurent series form.
The algebraically finite elements of $\unit(\mathfrak A[\mathsf z^{-1},\mathsf z])$ fall into various finiteness classes $F$
depending on the length of the elements $a_s$ or $(a_s)^{-1}$.
Let ${\mathfrak A}_F$ be the set of decompositions $\{a_j\}_{1\leq j\leq s}$ compatible with $F$.
Then Statement \ref{stat:linloop} can be augmented as follows:

\begin{state}\label{stat:subfinite}
For any finiteness class $F$, there is a smooth homotopy
\[K_F^{\mathrm e}:
{\mathfrak A}_F\times[0,1]\times[0,\pi/2]\rightarrow\unit(\mathcal K_{\mathbb Z}(\mathfrak A)[\mathsf z^{-1},\mathsf z]),\]
such that

(i) $K_F^{\mathrm e}(\tilde a,0,\theta)=K(a,\theta)$;

(ii) $K_F^{\mathrm e}(\tilde a,1,\theta)$ differs from $1_{\mathbb Z}$
in finitely many places (depending on $F$);

(iii) $K_F^{\mathrm e}(\tilde a,h,\pi/2)$ is constant in $h$;

(iv)
$K_F^{\mathrm e}(\tilde a,h,0)={\mathsf U}_{F}(\tilde a,h)\Lambda(\mathsf z,\mathsf Q)
{\mathsf U}_{F}(\tilde a,h)^{-1}\Lambda(\mathsf z,\mathsf Q)^{-1}.$
\newline
Here ${\mathsf U}_{F}(\tilde a,h)$ differs from $\mathsf U(a)$ in
a rapidly decreasing matrix.

In particular, $K_F^{\mathrm e}(\tilde a,1,0)$ yields a finite linearization of $a(\mathsf z)a(1)^{-1}$.
\end{state}

These statements are known, but in lesser generality,
in  various ways:
Statement \ref{stat:infstable}, as stated here in the smooth category (however, see \ref{stat:finite}),
follows from Cuntz, \cite{C2}, Section 2.
Statement \ref{stat:linloop} is
a quantitative version of the well-known linearization technique of Atiyah and Bott, \cite{AB};
but much resembling to the formulas of Pressley and Segal, \cite{PS}, Ch.~6, who work with
Hilbert-Schmidt matrices, instead of rapidly decreasing ones. Statement \ref{stat:toepcont}
comes from  the original Toeplitz  argument of Cuntz, \cite{C}, originally stated in the context of
$C^*$-algebras, but subsequently adapted to the smooth case, cf. also \cite{C2}.
One can also find some explicit homotopies in \cite{C3}.
Statement \ref{stat:finite} is useful,
because $*$-algebras are prominent in operator algebraic discussions; and the finite
perturbation category is the technically easiest setting to provide large
contractible spaces for the purposes of algebraic topology.
Statement \ref{stat:subfinite} amounts to an explicit computation in the less functorial but more concrete
setting of \cite{AB}.

The constructions presented here are improved versions of some
constructions which can be found in the author's thesis \cite{L}.
The author indebted to Prof. Richard B. Melrose, his advisor, for helpful discussions.
In fact, much of this content was motivated by the geometric idea of Melrose, Rochon \cite{MR}.
The author would also like to thank Prof. Joachim Cuntz, who called his attention to some related papers,
and Prof. Bal\'azs Csik\'os, for some useful advices.

\section{A general framework for computations }\label{sec:comp}

If $\mathfrak A$ is a not necessarily unital algebra, then one can consider
the semigroup $1+\mathfrak A$, with
elements of form $1+a$, $(a\in\mathfrak A)$, which multiply as $(1+a)(1+b)=1+(a+b+ab)$.
If $\mathfrak A$ is unital, then it is customary to identify $\mathfrak A$ and $1+\mathfrak A$
by the recipe
$a\in\mathfrak A\,\leftrightarrow\, 1-(1_{\mathfrak A}-a)\in1+\mathfrak A$.
This is also the situation
if there is a natural identity element which can be associated to $\mathfrak A$, like the
identity matrix in the case of matrix algebras.
The unit group $\unit(\mathfrak A)$ of $\mathfrak A$ is the
unit group  of the semigroup $1+\mathfrak A$, i.~e., it is the group of pairs
$(1+a,1+b)\in(1+\mathfrak A)\times(1+\mathfrak A)$
such that $(1+a)(1+b)=(1+b)(1+a)=1$; they multiply
$(1+a_1,1+b_1)(1+a_2,1+b_2)=((1+a_1)(1+a_2),(1+b_2)(1+b_1))$.
If $\mathfrak A$ is a topological ring, then the natural topology
on $\unit(\mathfrak A)$  comes from the product topology of
$(1+\mathfrak A)\times(1+\mathfrak A)$ by restriction.
As $1+a$ determines $1+b$, we write ``$1+a$'' instead of ``$(1+a,1+b)$''.
If $\phi:\mathfrak A\rightarrow\mathfrak B$
is a homomorphism, then it induces a homomorphism
$\unit\phi:\unit(\mathfrak A)\rightarrow\unit(\mathfrak B)$ defined by
$1+a\mapsto1+\phi(a)$.
We will write $\phi$ instead of $\unit\phi$.

In what follows, a ``locally convex vector space $\mathfrak A$''  means a
sequentially complete, Hausdorff, locally convex vector space $\mathfrak A$.
The completeness is essential for analytic purposes.
If the topology of $\mathfrak A$ is induced by a set $\Pi_{\mathfrak A}$ of seminorms,
then we assume that any positive integral combination of these seminorms also belongs
to the generating seminorm set.
A locally convex algebra $\mathfrak A$ is a locally convex vector space
with continuous bilinear multiplication. So, for each seminorm  $p\in\Pi_\mathfrak A$
there is an other seminorm $\tilde p\in\Pi_\mathfrak A$ such that for all $X_1,X_2\in\mathfrak A$
the inequality $p(X_1X_2)\leq\tilde p(X_1)\tilde p(X_2)$ holds.
An inductive locally convex vector space $\mathfrak A$ is an indexed family of
locally convex vector spaces $\{\mathfrak A_\lambda\}_{\lambda\in\Lambda}$
such that the following holds:
$\Lambda$ is an upward directed partially ordered set, i. e.
for all $\lambda,\mu\in\Lambda$ there is an element $\nu\geq\lambda,\mu$.
For all $\mu\geq\lambda$ there exist continuous inclusions
$T_\mu^\lambda:\mathfrak A_\lambda\rightarrow\mathfrak A_\mu$; and for
$\nu\geq\mu\geq\lambda$ one has $T_\nu^\mu\circ T_\mu^\lambda=T_\nu^\lambda$.
Now, $\mathfrak A$ is an inductive locally convex algebra if for each $\lambda,\mu\in\Lambda$ there is
an element $\mathrm{prod}(\lambda,\mu)\in\Lambda$, and for
$\nu\geq\mathrm{prod}(\lambda,\mu)$  bilinear products
$M_{\lambda,\mu}^\nu:\mathfrak A_\lambda\times\mathfrak A_\mu\rightarrow\mathfrak A_{\nu}$
compatible with the inclusions and the usual algebraic prescriptions are given.
An element of $\mathfrak A$ is an element of $\bigcup_{\lambda\in\Lambda}\mathfrak A_\lambda$ making
identifications along the inclusion maps.
Then $\mathfrak A$ will be an algebra.
endowed with an ``inductive'' topology
coming from the filtration $\{\mathfrak A_\lambda\}_{\lambda\in\Lambda}$, such that the
vector space structure respects the filtration but the algebra structure does not.
If  the spaces $\mathfrak A$ and $\mathfrak B$ have inductive topologies with filtrations
$\{\mathfrak A_\lambda\}_{\lambda\in\Lambda}$
and $\{\mathfrak B_\mu\}_{\mu\in M}$, then a map $\phi:\mathfrak A\rightarrow\mathfrak B$
is continuous if for each $\lambda\in\Lambda$ there is an element $\mu\in M$
such that there is a continuous map $\phi_\lambda:\mathfrak A_\lambda\rightarrow
\mathfrak B_\mu$, which is set-theoretically a restriction of $\phi$.

Suppose that $\Theta_1,\Theta_2$ are sets and $\mathfrak V$ is a vector space.
Then a $\mathfrak V$-valued $\Theta_1$ times $\Theta_2$ matrix is just a formal sum
${s=\sum_{a\in\Theta_1, b\in\Theta_2}s_{a,b}\mathbf e_{a,b}\in
\mathcal M_{\Theta_1,\Theta_2}(\mathfrak V)}$
with coefficients $s_{a,b}$ from $\mathfrak V$.
We write $\mathcal
M_\Theta(\mathfrak V)$ instead of $\mathcal M_{\Theta,\Theta}(\mathfrak V)$ and use similar notation
for other spaces as well.
For column and row matrices, we use the notation $\mathbf e_a=\mathbf e_{a,*}$ and
$\mathbf e_b^\top=\mathbf e_{*,b}$ respectively, and we make the formal identification
$\mathbf e_{a,b}=\mathbf e_a\otimes \mathbf e_b^\top$. For
column spaces, we use the notation
$\mathcal S(\Theta;\mathfrak V)=\mathcal M_{\Theta,\{*\}}(\mathfrak V)$.
We use the notation
${1_\Theta=\sum_{\theta\in\Theta}\mathbf e_{\theta,\theta}}$,
and in general circumstances we consider the identity matrix $1_\Theta$ as the
adjoint unit in any non-unital $\Theta$ times $\Theta$ matrix algebra.
If $s_i\in\mathfrak V$, $i\in\mathbb Z$ are given then
\begin{equation}\Diag(\ldots s_{-2},s_{-1}|s_0,s_1,s_2,\ldots)=
\sum_{i\in\mathbb Z}s_i\mathbf e_{i,i}\in
\mathcal M_{\mathbb Z}(\mathfrak V)
\label{eq:diag}
\end{equation}
is the corresponding diagonal matrix.
$\Diag(s_0,s_1,s_2,\ldots)\in\mathcal M_{\mathbb N}(\mathfrak V)$, similarly.
For $a\in\mathfrak A$, we define the matrices
$\mathsf E_{\mathbb N}(a)=a\mathbf e_{00}\in\mathcal K_{\mathbb N}(\mathfrak A)$
and
$\mathsf E_{\mathbb Z}(a)=a\mathbf e_{00}\in\mathcal K_{\mathbb Z}(\mathfrak A)$.
Then, as usual, for $\tilde a=1+a\in1+\mathfrak A$, we extend these maps as
$\mathsf E_{\mathbb N}(\tilde a)=1_{\mathbb N}+\mathsf E_{\mathbb N}(a)$
and
$\mathsf E_{\mathbb Z}(\tilde a)=1_{\mathbb Z}+\mathsf E_{\mathbb Z}(a)$;
i. e., for $\tilde a\in1+\mathfrak A$, it yields
$\mathsf E_{\mathbb N}(\tilde a)=\Diag(\tilde a,1,1,\ldots)$,
and
$\mathsf E_{\mathbb Z}(\tilde a)=\Diag(\ldots,1,1\,|\,\tilde a,1,1,\ldots)$.

On the set $\mathbb N$ of natural numbers, there is the natural space
$\mathcal S^{\infty}(\mathbb N;\mathbb R)^\star$, i.~e. the
space of multiplicatively invertible polynomially growing functions.
A countable set $\Theta$ is called a set of polynomial growth if it is endowed with
a set of
functions $\mathcal S^{\infty}(\Theta;\mathbb R)^\star$
from $\Theta$ to $\mathbb R$ such that there is a bijection
$\omega:\Theta\rightarrow\mathbb N$  such that
 $\omega^*\mathcal S^{\infty}(\mathbb N;\mathbb R)^\star=
\mathcal S^{\infty}(\Theta;\mathbb R)^\star$.
It is notable that $\mathbb N\times\mathbb N$ and $\mathbb N\,\dot\cup\,\mathbb N$ are
sets of polynomial growth naturally;
and that way we can define the direct product $\Theta_1\times\Theta_2$ and direct
sums $\Theta_1\,\dot\cup\,\Theta_2$ of sets of
polynomial growth $\Theta_1$ and $\Theta_2$.
In what follows, the sets of polynomial growth we use will be like
$\mathbb N, \mathbb Z,$ or
$\{1,\ldots,n\}\times \mathbb Z$, where the description of the relevant function
spaces is evident, so it will not be detailed.
The main point is that  a set $\Theta$ of polynomial growth is just
like $\mathbb N$ for practical purposes.
If $\Theta_1,\Theta_2$ are sets of polynomial growth, and $\mathfrak V$ is a locally convex
vector space, then we can define some matrix spaces as follows:

(a)
With functions $F:\Pi_{\mathfrak V}\rightarrow  \mathcal S^{\infty}(\Theta_1;\mathbb R)^\star\times\mathcal S^{\infty}(\Theta_2;\mathbb R)^\star$, the filtering spaces
\begin{multline}
\mathcal M^{\infty,\infty}_{\Theta_1,\Theta_2}(\mathfrak V)_F=
\biggl\{ s\in \mathcal M_{\Theta_1,\Theta_2}(\mathfrak V)
\,:\,\forall p\in\Pi_{\mathfrak V}
\text{\qquad}\\ |s|_{\frac1{F_1(p)},p,\frac1{F_2(p)}}=\sum_{(a,b)\in\Theta_1\times\Theta_2}
\left|\frac1{F_1(p)(a)}\right|\,p(s_{a,b})\left|\frac1{F_2(p)(b)}\right|<+\infty\biggr\}\notag
\end{multline}
form the inductive locally  convex space
 $\mathcal M^{\infty,\infty}_{\Theta_1,\Theta_2}(\mathfrak V)$.

(b)
With functions $F:\Pi_{\mathfrak V}\times\mathcal S^{\infty}(\Theta_2;\mathbb R)^\star
 \rightarrow \mathcal S^{\infty}(\Theta_1;\mathbb R)^\star$,
 the filtering spaces
\begin{multline}
\mathcal M^{\infty,-\infty}_{\Theta_1,\Theta_2}(\mathfrak V)_F=
\biggl\{ s\in \mathcal M_{\Theta_1,\Theta_2}(\mathfrak V)
\,:\,\forall p\in\Pi_{\mathfrak V} \forall g\in \mathcal S^{\infty}(\Theta_2;\mathbb R)^\star
\text{\qquad}\\ |s|_{\frac1{F(p,g)},p,g}=\sum_{(a,b)\in\Theta_1\times\Theta_2}
\left|\frac1{F(p,g)(a)}\right|\,p(s_{a,b})|g(b)|<+\infty\biggr\}\notag
\end{multline}
form the inductive locally  convex  space
$\mathcal M^{\infty,-\infty}_{\Theta_1,\Theta_2}(\mathfrak V)$.
We can define the space
$\mathcal M^{-\infty,\infty}_{\Theta_1,\Theta_2}(\mathfrak V)$ similarly.

(c) We define
\begin{multline}
\mathcal M^{-\infty,-\infty}(\Theta_1,\Theta_2;\mathfrak V)=
\biggl\{ s\in \mathcal M(\Theta_1,\Theta_2;\mathfrak V)
\,:\,\forall p\in\Pi_{\mathfrak V} \forall f\in \mathcal S^{\infty}(\Theta_1;\mathbb R)^\star
\\ \forall g\in \mathcal S^{\infty}(\Theta_2;\mathbb R)^\star
\text{\qquad }|s|_{f,p,g}=\sum_{(a,b)\in\Theta_1\times\Theta_2}
|f(a)|\,p(s_{a,b})|g(b)|<+\infty\biggr\}.\notag
\end{multline}

(d) It is natural to define
$\Psi_{\Theta_1,\Theta_2}(\mathfrak V)=
\mathcal M^{-\infty,\infty}_{\Theta_1,\Theta_2}(\mathfrak V)
\cap \mathcal M^{\infty,-\infty}_{\Theta_1,\Theta_2}(\mathfrak V),$
the space of matrices of ``pseudodifferential size''.

If $\mathfrak A\times\mathfrak B\rightarrow \mathfrak C$ is a continuous bilinear pairing
between locally convex spaces, then we have induced continuous pairings
$\mathcal M^{\mathsf X,\infty}_{\Theta_1,\Theta_2}(\mathfrak A)
\times \mathcal M^{-\infty,\mathsf Y}_{\Theta_2,\Theta_3}(\mathfrak B)
\rightarrow \mathcal M^{\mathsf X,\mathsf Y}_{\Theta_1,\Theta_3}(\mathfrak C)$,
$\Psi_{\Theta_1,\Theta_2}(\mathfrak A)
\times  \Psi_{\Theta_2,\Theta_3}(\mathfrak B)
\rightarrow \Psi_{\Theta_1,\Theta_3}(\mathfrak C),$
etc.
So come the algebra and module structures associated to matrices.
Instead of $\mathcal M^{-\infty,-\infty}_{\Theta_1,\Theta_2}(\mathfrak A)$
we will use the shorter notations $\mathcal K_{\Theta_1,\Theta_2}(\mathfrak A)$.
Instead of $\mathcal M^{\pm\infty,\mathsf X}_{\Theta,\{*\}}(\mathfrak A)$ it is reasonable to use
$\mathcal S^{\pm\infty}(\Theta;\mathfrak A)$, which is a consistent extension of our earlier notation.
There are natural isomorphisms like
$\mathcal K_{\Theta_1\times\Theta_1',\Theta_2\times\Theta_2'}(\mathfrak A)
\simeq\mathcal K_{\Theta_1,\Theta_2}(\mathcal K_{\Theta_1',\Theta_2'}(\mathfrak A))$, etc.
One often uses is relabeling of matrices, which is as follows:
Suppose that $\omega:\Omega\rightarrow\Omega'$
is a map between sets of polynomial growth, such that
$\omega^*\mathcal S^{\infty}(\Omega',\mathbb R)^\star=\mathcal S^{\infty}(\Omega,\mathbb R)^\star$.
This includes the case when $\omega$
is an isomorphism of sets of polynomial growth,
and also the natural inclusions $\iota:\Omega\rightarrow\Omega'\dot\cup\Omega''$, where $\Omega''$ is
finite or an other set of polynomial growth.
Let us now consider the matrix
$R_\omega=\sum_{\alpha\in\Omega}\mathbf e_{\alpha,\omega(\alpha)}\in\Psi_{\Omega,\Omega'}(\mathbb R)$.
Then, for a matrix $A\in\mathcal M^{\mathsf X,\mathsf Y}_\Theta(\mathfrak A)$ or
$\Psi_\Theta(\mathfrak A)$, we can take the matrix
$r_\omega(A)=R_{\omega}^\top AR_\omega$ which is a matrix of  the same kind as $A$
but $\Omega$ is replaced by $\Omega'$.
This relabeling $r_\omega$ is a continuous, smooth operation,
which is an isomorphism if $\omega$ is an isomorphism.

The advantage of the spaces $\mathcal M^{\mathsf X,\mathsf Y}_{\Theta_1,\Theta_2}(\mathfrak V)$
is that they are sufficiently large for the purposes of arithmetic calculations.
In what follows,  only the algebras $\mathcal K$
will be used explicitly.  On the other hand, all computations,
except in Section \ref{sec:otheralg}
will be governed by the principle
every matrix expression will be understood as an element of $\Psi_{\Omega_1,\Omega_2}(\mathfrak A)$,
where $\Omega_i$ are sets of polynomial growth, and $\mathfrak A$ is a locally convex algebra;
but we always hope that  our expressions will yield results
which turn out to be continuous in stronger topologies.

\section{The environment of cyclic and Toeplitz algebras}\label{sec:CT}
\subsection*{Cyclic and Toeplitz algebras}
In what follows, let $\overline{\mathbb N}=\mathbb Z\setminus\mathbb N$, so
$\mathbb Z=\overline{\mathbb N}\,\dot\cup\,\mathbb N$.
We make a canonical correspondence between $\mathbb N$ and
$\overline{\mathbb N}$ by relabeling every $n$ to $-1-n$.
We can consider every $\mathbb Z\times\mathbb Z$ matrix $U$
as a $2\times 2$ matrix of $\mathbb N\times\mathbb N$ matrices:
\[U=\left[\begin{array}{c|c}
U|_{\overline{\mathbb N}\times\overline{\mathbb N}} &
U|_{\overline{\mathbb N}\times{\mathbb N}} \\\hline
U|_{{\mathbb N}\times\overline{\mathbb N}} &
U|_{{\mathbb N}\times{\mathbb N}}
\end{array}\right]\simeq
\begin{bmatrix}U^{--}&U^{-+}\\U^{+-}&U^{++}\end{bmatrix},\]
such that  the matrix entries on the right side are $\mathbb N\times\mathbb N$ matrices
obtained by the correspondence explained above.

An element $a=\sum_{i\in\mathbb Z}a_i\mathbf e_i\in\mathcal
S^{-\infty}(\mathbb Z;\mathfrak A)$ can and will, in general,  be identified
with the Laurent series $\sum_{i\in\mathbb Z}a_i\mathsf z^i\in\mathfrak A[\mathsf z^{-1},\mathsf z]$ with rapidly decreasing coefficients.
We call this algebra  the
algebra of cyclic loops, in contrast to the algebra of proper loops
$\mathcal C^{\infty}(\mathrm S^1;\mathfrak A)$.
Elements $a=\sum_{i\in\mathbb Z}a_i\mathsf z^i\in\mathfrak A[\mathsf z^{-1},\mathsf z]$ can be
represented by $\mathbb Z\times\mathbb Z$ matrices
\begin{equation}
\mathsf U(a)=\sum_{n,m\in\mathbb Z}a_{n-m}\mathbf e_{n,m}=
\left[\begin{array}{ccc|ccc}
\ddots&\ddots&\ddots&\ddots&\ddots&\ddots\\
\ddots&a_0&a_{-1}&a_{-2}&a_{-3}&\ddots\\
\ddots&a_1&a_0&a_{-1}&a_{-2}&\ddots\\\hline
\ddots&a_2&a_1&a_0&a_{-1}&\ddots\\
\ddots&a_3&a_2&a_1&a_0&\ddots\\
\ddots&\ddots&\ddots&\ddots&\ddots&\ddots
\end{array}\right].
\label{eq:u}
\end{equation}

If $\mathcal W_{\mathbb Z}(\mathfrak A)$ is the image set
of $\mathfrak A[\mathsf z^{-1},\mathsf z]$ under $\mathsf U$, then it is a subset
of $\Psi(\mathbb Z;\mathfrak V)$ algebraically, but we put the topology of
$\mathcal S^{-\infty}(\mathbb Z;\mathfrak V)$ to it. If $\mathfrak A$ is a locally
convex algebra, then $\mathsf U:\mathfrak A[\mathsf z^{-1},\mathsf z]\rightarrow\mathcal W_{\mathbb Z}(\mathfrak A)$
is an isomorphism of algebras.
When it comes to the $2\times 2$ decomposition as explained above,
in order to simplify the notation, we will just write $\mathsf W(a)$ instead
of $\mathsf U(a)^{++}$, and $\mathsf Y(a)$ instead of $\mathsf U(a)^{+-}$.
If $a=a(\mathsf z)$ then, with some abuse of notation, we also write $a^\top=a(\mathsf z^{-1})$.
Then
\[\mathsf U(a)=\begin{bmatrix}
 \mathsf W(a^\top)& \mathsf Y(a^\top)\\
 \mathsf Y(a) & \mathsf W(a)
\end{bmatrix}.\]
So $\mathsf W(a)$ is the infinite Toeplitz matrix
associated to $a$, and $\mathsf Y(a)$ is the infinite Hankel matrix
associated to ("the positive part" of) $a$.

As far as the linear structure is concerned, we could have just used
the matrices $\mathsf W(a)$ to represent the elements $a$. The difference is that
in terms of matrix multiplication
$\mathsf W(a)\mathsf W(b)=\mathsf W(ab)-\mathsf Y(a)\mathsf Y(b^\top),$
so there is an ``anomalous'' term
$-\mathsf Y(a)\mathsf Y(b^\top)\in\mathcal K_{\mathbb N}(\mathfrak A )$.
One can see that
algebraically $\mathcal W_{\mathbb N}(\mathfrak A )
\cap\mathcal K_{\mathbb N}(\mathfrak A )=0$.
Hence it is reasonable to define the Toeplitz algebra
\[\mathcal T_{\mathbb N}(\mathfrak A )=
\mathcal W_{\mathbb N}(\mathfrak A )+
\mathcal K_{\mathbb N}(\mathfrak A ),\]
which is topologically just
$\mathcal W_{\mathbb N}(\mathfrak A )\oplus
\mathcal K_{\mathbb N}(\mathfrak A )$ but with the algebraic  product rule
$(\mathsf W(a)+p)((\mathsf W(b)+q))=\mathsf W(ab)+(
-\mathsf Y(a)\mathsf Y(b^\top)+\mathsf W(a)q+p\mathsf W(b)+pq),$
induced from the matrix structure. Algebraically,
$\mathcal T_{\mathbb N}(\mathfrak A )$
is just a subset of $\Psi(\mathbb N;\mathfrak A )$ but a locally convex algebra.
So, one can see that there is a short exact sequence of algebras
$0\rightarrow \mathcal K_{\mathbb N}(\mathfrak A )
\xrightarrow\iota \mathcal T_{\mathbb N}(\mathfrak A )
\xrightarrow\sigma \mathcal W_{\mathbb Z}(\mathfrak A )\rightarrow 0$.
The map $\iota$ is the inclusion of the ideal of rapidly decreasing matrices
into the Toeplitz algebra,
while $\sigma$ is the symbol map.
In what follows, we rather consider the value of the symbol map as an element of
$\mathfrak A[\mathsf z^{-1},\mathsf z]$,
so we have the symbol homomorphism
\[\sigma: \mathcal T_{\mathbb N}(\mathfrak A )\rightarrow
\mathfrak A[\mathsf z^{-1},\mathsf z].\]
We can naturally extend this symbol map to unit groups as we have seen.

For technical reasons, we  define the algebra
\[\mathcal T_{\mathbb Z}(\mathfrak A )=
\begin{bmatrix}
\mathcal T_{\mathbb N}(\mathfrak A )& \mathcal K_{\mathbb N}(\mathfrak A )\\
\mathcal K_{\mathbb N}(\mathfrak A ) & \mathcal T_{\mathbb N}(\mathfrak A )
\end{bmatrix},\]
which is also naturally a locally convex algebra. Then
$\mathcal W_{\mathbb Z}(\mathfrak A)\subset \mathcal W_{\mathbb Z}(\mathfrak A )+
\mathcal K_{\mathbb Z}(\mathfrak A )\subset \mathcal T_{\mathbb Z}(\mathfrak A )$.
For the sake of notational convenience, we define the block matrix
\[\widehat{\mathsf  U}(a)=\left[\begin{array}{c|cc}
\mathsf W(a^\top)&&-\mathsf Y(a^\top)
\\\hline
&0&\\
-\mathsf Y(a)&&\mathsf W(a)
\end{array}\right]\in \mathcal T_{\mathbb Z}(\mathfrak A ).\]
We remark that for $\tilde a\in\mathtt U(\mathfrak A[\mathsf z^{-1},\mathsf z])$ an ``$1$'' appears in the place of  ``$0$''.

Elements of $\mathcal T_{\mathbb Z}(\mathfrak A)$ have
two symbols; one belonging to the lower right quadrant,
and one belonging to the upper left quadrant.
It is a small but important observation regarding $
\mathsf U(a)\in\mathcal T_{\mathbb Z}(\mathfrak A)$
that the Toeplitz element in the lower right quadrant has
symbol $a=a(\mathsf z)$, but the Toeplitz element in the upper left quadrant has symbol
$a^\top=a(\mathsf z^{-1})$.

One can also see that there are natural isomorphisms like
$\mathcal T_{\mathbb N}(\mathcal K_\Omega(\mathfrak A))\simeq
\mathcal K_\Omega(\mathcal T_{\mathbb N}(\mathfrak A)),$
etc.
In fact, all of our matrix space constructions considered as functors  are naturally ``commutative''.

Let $\mathfrak A[\mathsf z^{-1},\mathsf z]^{\mathrm{po}}$ be the set of pointed loops, i. e., where $a(1)=0$.
Then the elements $\tilde a\in\unit(\mathfrak A[\mathsf z^{-1},\mathsf z]^{\mathrm{po}}$)
are those for which $\tilde a(1)=1$.
These pointed spaces are closed subspaces of the unpointed spaces.
We can define the pointed Toeplitz algebra
$\mathcal T_{\mathbb N}(\mathfrak A)^{\mathrm{po}}$ similarly,
 the symbols are pointed there.
\subsection*{The Bott involution map}
In what follows, we use the abbreviation
$\Lambda(a,b)=\frac12(1+a+b-ab)$.
Let
$\mathsf Q=\begin{bmatrix}-1_{{\mathbb N} }&\\&1_{{\mathbb N} }\end{bmatrix}$.
 Then
$\Lambda(\mathsf z,\mathsf Q)=
\begin{bmatrix}\mathsf z1_{{\mathbb N} }&\\&1_{{\mathbb N}}\end{bmatrix}.$
We also use the delta-function $\delta_{n,m}$, which is 1 if $n=m$, and it is $0$ otherwise.

If $a\in\unit(\mathfrak A[\mathsf z^{-1},\mathsf z])$, then we define the ``Bott'' involution
\[\mathsf B(a)=\mathsf U(a)\mathsf Q\mathsf U(a)^{-1}\in\mathsf Q+
\mathcal K_{\mathbb Z}(\mathfrak A)\]
(cf. the symbols).

\subsection*{``Shifting rotations''}
Our natural deformation parameter variable, in general, will be $\theta\in[0,\pi/2]$,
or, more generally,  $\theta\in\mathrm S^1=\mathbb R/2\pi\mathbb Z$.
In order to save space, we often use $t=\sin \theta$ and $s=\cos \theta$ instead.
It is useful to keep in mind that $s^2=1-t^2$.
For $\theta\in\mathrm S^1$,  we define the matrices
\[\mathsf C(\theta)=\begin{bmatrix}
s&ts&t^2s&t^3s&\cdots\\
-t & s^2& ts^2 & t^2s^2  &\ddots\\
&-t & s^2 & ts^2 & \ddots\\
&&-t & s^2 & \ddots\\
&&&-t & \ddots\\
&&&&\ddots\\
\end{bmatrix}\in
\mathcal T_{\mathbb N}(\mathbb R).\]

\begin{lemma}\label{lem:artkey}
Let $\mathsf C(\theta)^\dag$ denote the transpose of $\mathsf C(\theta)$. Then

(a)\[\mathsf C(\theta)^\dag\mathsf C(\theta)=1_{\mathbb N}.\]

(b)
\[\mathsf C(\theta)\mathsf C(\theta)^\dag=
-\delta_{t,1}\mathbf e_{0,0}-\delta_{t,-1}\mathbf e_{0,0}+1_{\mathbb N}.\]

(c)
\[\mathsf C(\theta)\mathbf e_{n,m}\mathsf C(\theta)^\dag=
\begin{bmatrix}
t^{n+m}s^2&t^{n+m-1}s^3&\cdots& t^{n}s^3
& -t^{n+1}s &\\
t^{n+m-1}s^3&t^{n+m-2}s^4&\cdots& t^{n-1}s^4
& -t^{n}s^2 &\\
\vdots &\vdots &\ddots
& \vdots &\vdots\\
t^{m}s^3&t^{m-1}s^4&\cdots&s^4
& -ts^2 &\\
 -t^{m+1}s&-t^{m}s^2&\cdots& -ts^2&t^2\\
 &&&&&0\\&&&&&&\ddots
\end{bmatrix}.\]

(d) For $n>0$,
\[\mathsf C(\theta)\mathsf W(\mathsf z^n)\mathsf C(\theta)^\dag=
-\delta_{t,1}\mathbf e_{0,0}-(-1)^n\delta_{t,-1}\mathbf e_{0,0}+\begin{bmatrix}
t^n&\\
t^{n-1}s&\\
\vdots\\
ts&\\
s&\\
&1\\
&&1\\
&&&\ddots
\end{bmatrix},\]
and
$\mathsf C(\theta)\mathsf W(\mathsf z^{-n})\mathsf C(\theta)^\dag$
is the transpose of the matrix above.
\begin{proof}
Direct computation.
\end{proof}
\end{lemma}
\begin{remark}
The presence of the terms $\delta_{t,1},\delta_{t,-1}$ might be surprising at first sight.
It reflects the phenomenon that in a topological algebra one cannot simultaneously topologize
the families  $\mathsf C(\theta)$ and $\mathsf C(\theta)^\dag$ correctly.
In fact, the $\mathsf C(\theta)$'s are isometries for $-1<t<1$, but they are just partial isometries for $t=\pm1$.
\end{remark}

\section{Stabilizing homotopies}\label{sec:stabilize}
\begin{proposition}\label{prop:rapid} The continuous map
\[T_{\mathcal K}:\mathcal K_{\mathbb N}(\mathfrak A)\times \mathrm S^1
\rightarrow\mathcal K_{\mathbb N}(\mathfrak A)\]
given by
\[A,\theta\mapsto T_{\mathcal K}(A,\theta)=\mathsf C(\theta)A\mathsf C(\theta)^\dag\]
is smooth in $\theta$.
It yields a family of endomorphisms of $\mathcal K_{\mathbb N}(\mathfrak A)$ when $\theta$ is fixed.
These are isomorphisms for $-1<t<1$, and  closed injective endomorphisms for $t=\pm1$.
In particular, for $\theta=0$ $(t=0)$,
\[T_{\mathcal K}(A,0)=A=\begin{bmatrix}
a_{11} & a_{12} &\cdots\\
a_{21} & a_{22} &\cdots\\
\vdots&\vdots&\ddots
\end{bmatrix};\]
but for $\theta=\pm\pi/2$ ($t=\pm 1$),
\[T_{\mathcal K}(A,\pm\pi/2)=
\begin{bmatrix}0\\
&a_{11} & a_{12} &\cdots\\
&a_{21} & a_{22} &\cdots\\
&\vdots&\vdots&\ddots
\end{bmatrix}.\]
\begin{proof}
Well-definedness and smoothness follows from Lemma \ref{lem:artkey}.c.
Lem\-ma \ref{lem:artkey}.a implies that we have a family of endomorphisms.
Furthermore, it also shows that
$A=\mathsf C(\theta)^\dag T_{\mathcal K}(A,\theta)\mathsf C(\theta)$;
from which the statement about the nature of the endomorphisms follows easily.
\end{proof}
\end{proposition}
Hence, taking $\theta\in[0,\pi/2]$, we see that the deformation $T_{\mathcal K}$ does
indeed realize a stabilizing homotopy, even if only with one ``extra dimension''.
Nevertheless, after this, stabilization becomes a matter of standard tricks:
\begin{cor}[$\Rightarrow$ Statement \ref{stat:infstable}]  \label{cor:infstable}
Let $\Omega_1$ and $\Omega_2$ be sets of polynomial growth;
$\Omega=\Omega_1\,\dot\cup\,\Omega_2$, and let $\omega:\Omega\rightarrow\Omega_1\subset\Omega$ be the composition
of an isomorphism $\Omega\simeq \Omega_1$ and the natural inclusion
$\Omega_1\rightarrow\Omega_1\,\dot\cup\,\Omega_2$. Then we claim:

There is a smooth map
$\widehat T_{\mathcal K}:\mathcal K_\Omega(\mathfrak S)\times[0,\pi/2]
\rightarrow\mathcal K_\Omega(\mathfrak S)$
such that it yields a family of endomorphisms of $\mathcal K_\Omega(\mathfrak S)$,
which are isomorphisms for $\theta\in[0,\pi/2)$, and a
closed injective endomorphism for $\theta=\pi/2$, such that
$\widehat T_{\mathcal K}(A,0)=\id_{\mathcal K_\Omega(\mathfrak S)}$ and
$\widehat T_{\mathcal K}(A,\pi/2)=r_\omega$.
The map $\widehat T_{\mathcal K}$ extends to unit groups naturally.
\begin{proof}
Take $\mathfrak A=\mathcal K_{\mathbb N}(\mathfrak S)$ in the previous statement.
It yields our statement with $\Omega_1=(\mathbb N\setminus\{0\})\times\mathbb N$,
$\Omega_2=\{0\}\times\mathbb N$,
$\omega((n,m))=(n+1,m)$.
Now, using an appropriate relabeling $r_\eta$ of $\Omega$ we obtain the general statement.
\end{proof}
\end{cor}
\begin{remark}
Another way to achieve stabilization by many dimensions is to ``quantize'' $\mathsf C(\theta)$,
see \cite{LL}.
\end{remark}

Due to the multiplicative structure, the concatenation of group valued homotopies
is particularly simple: If $f,g:Y\times[0,1]\rightarrow G$ yield homotopies $f_0\simeq f_1$, $g_0\simeq g_1$
where $f_1=g_0$, then $h(y,t)=f(y,t)f(y,1)^{-1}g(y,t)$ yields a homotopy between $f_0$ and $g_1$.
Then polynomial/smooth homotopies yield polynomial/smooth
homotopies, and the operation is associative; in contrast to concatenation by reparametrization.
Using this observation and the stabilizing homotopies above, one can easily prove
\begin{cor}\label{cor:stablehom}
Let $\Omega_1,\Omega_2$ be sets of polynomial growth, and let
$\iota_1:\Omega_1\rightarrow\Omega_1\,\dot\cup\,\Omega_2$
be the natural inclusion.
Assume that
$H:X\times [0,1]\rightarrow\unit(\mathcal K_{\Omega_1\,\dot\cup\,\Omega_2}(\mathfrak S))$
is a smooth homotopy with maps $f_0, f_1:X\rightarrow
\unit(\mathcal K_{\Omega_1}(\mathfrak S))$
such that
$H_0=r_{\iota_1}(f_0)$ and $H_1=r_{\iota_1}(f_1)$.
Then we claim that there is a smooth homotopy
$f:X\times [0,1]\rightarrow\unit(\mathcal K_{\Omega_1}(\mathfrak S))$
between $f_0$ and $f_1$. This $f$ can be chosen so that there is a
smooth homotopy between $H$ and $r_{\iota_1}(f)$ relative to endpoints.
In other words: ``In stable algebras stable homotopies can be reduced to ordinary homotopies.''
\qed
\end{cor}

\section{Linearization of cyclic loops}
\begin{point}
Let $\mathsf v$ be a cyclic formal variable, and take
$\mathsf V=\sum_{n\in\mathbb N}\mathsf v^{n}\mathbf e_{n,n}$.
Furthermore, take
$\displaystyle{\mathsf G(\theta,\mathsf v)=
\begin{bmatrix}1_{\mathbb N}&\\&\mathsf V^{-1}\mathsf C(\theta)\mathsf V\end{bmatrix}}
$ and $\displaystyle{
\mathsf G^\dag(\theta,\mathsf v)=
\begin{bmatrix}1_{\mathbb N}&\\&\mathsf V^{-1}\mathsf C(\theta)^\dag\mathsf V\end{bmatrix}
.}$
For $a\in\mathfrak A[\mathsf z^{-1},\mathsf z]$, we define
\[U(a,\theta,\mathsf v)=
\delta_{t,1}a(\mathsf v)\mathbf e_{00}+\delta_{t,-1}a(-\mathsf v)\mathbf e_{00}
+\mathsf G(\theta,\mathsf v)\mathsf U(a)\mathsf G^\dag(\theta,\mathsf v).\]
\begin{point}\label{po:udesc}
For $n>0$, this definition yields
\begin{multline}
U(\mathsf z^n, \theta,\mathsf v)=\\=
\left[\begin{array}{ccccccc|ccc}
\ddots&&&&&&&&&\\
&1&&&&&&&&\\
\hline
&&s&ts\mathsf v&t^2s\mathsf v^2&\cdots&t^{n-1}s\mathsf v^{n-1}&t^n\mathsf v^n&&\\
&&-t\mathsf v^{-1}&s^2&ts^2\mathsf v&\cdots&t^{n-2}s^2
           \mathsf v^{n-2}&t^{n-1}s\mathsf v^{n-1}&&\\
&&&\ddots&\ddots&\ddots&\vdots&\vdots&&\\
&&&&\ddots&\ddots&ts^2\mathsf v&t^2s\mathsf v^2&&\\
&&&&&\ddots&s^2&ts\mathsf v&&\\
&&&&&&-t\mathsf v^{-1}&s&&\\
&&&&&&&&1&\\
&&&&&&&&&\ddots
\end{array}\right]
\notag\end{multline}
$U(1, \theta,\mathsf v)=1_{\mathbb Z},$ and
$U(\mathsf z^{-n}, \theta,\mathsf v)=U(\mathsf z^{n}, \theta,\mathsf v^{-1})^\top$;
i. e. $U(\mathsf z^{n}, \theta,\mathsf v)$ is just a rather nice
perturbation / deformation of $\mathsf U(\mathsf z^n)$.
\end{point}
\begin{lemma}\label{lem:linearkey}
The continuous map
\[U:\mathfrak A[\mathsf z^{-1},\mathsf z]\times\mathrm S^1\rightarrow
\mathcal W_{\mathbb Z}(\mathfrak A)+
\mathcal K_{\mathbb Z}(\mathfrak A)[\mathsf v^{-1},\mathsf v]\]
defined by
\[ a,\theta\mapsto U(a,\theta,\mathsf v)\]
is smooth in $\theta$. It yields a family of homomorphisms with fixed $\theta$.
The symbols remain constant.
For $\theta=0$,
\[U(a,0,\mathsf v)=
\begin{bmatrix}
\mathsf W(a^\top) & \mathsf Y(a^\top)\\
\mathsf Y(a) & \mathsf W(a)
\end{bmatrix}
=\mathsf U(a);\]
while for $\theta=\pi/2$,
\[U(a,\pi/2,\mathsf v)=\left[\begin{array}{c|cc}
\mathsf W(a^\top)&&-\mathsf v\mathsf Y(a^\top)
\\\hline
&a(\mathsf v)&\\
-\mathsf v^{-1}\mathsf Y(a)&&\mathsf W(a)
\end{array}\right].\]
\begin{proof}
This is immediate from \ref{po:udesc} by taking linear combinations.
\end{proof}
\end{lemma}
 Considering  $a\in\mathtt U(\mathfrak A)$, and the natural extension to the unit group,
$U(a,\pi/2,\mathsf v)=\mathsf E_{\mathbb Z}(a(\mathsf v))\Lambda(\mathsf v,\mathsf Q)
\widehat{\mathsf U}(a)\Lambda(\mathsf v,\mathsf Q)^{-1}$ can be written.

\begin{proposition}[$\Rightarrow$ Statement \ref{stat:linloop}]\label{prop:linear}
The continuous map
\[K:\unit(\mathfrak A[\mathsf z^{-1},\mathsf z])\times\mathrm S^1\rightarrow
\unit(\mathcal K_{\mathbb Z}(\mathfrak A)[\mathsf v^{-1},\mathsf v]^{\pointed})\]
defined by
\[ a,\theta\mapsto K(a,\theta,\mathsf v)=
U(a,\theta,\mathsf v)\Lambda(\mathsf v,\mathsf Q)
U(a,\theta,1)^{-1}\Lambda(\mathsf v,\mathsf Q)^{-1}\]
is smooth in the variable $\theta$.
Here
\[K(a,0,\mathsf v)=
\Lambda(\mathsf v,\mathsf B(a))\Lambda(\mathsf v,\mathsf Q)^{-1},\qquad
K(a,\pi/2,\mathsf v)=
\mathsf E_{\mathbb Z}(a(\mathsf v)a(1)^{-1}).\]
\begin{proof} The statement follows immediately from the previous lemma.\end{proof}
\end{proposition}
\end{point}
\begin{remark}
When it comes to the linearization of not pointed loops but the
``cocycle'' $a(\mathsf z)a(\mathsf w)^{-1}$, then one can use the linearizing ``cocycle''
$K^{\mathrm c}(a,\theta,\mathsf z,\mathsf w)=U(a,\theta,\mathsf z)\Lambda(\mathsf z\mathsf w^{-1},
\mathsf Q)$ $U(a,\theta,\mathsf w)^{-1}$.
It yields
$K^{\mathrm c}(a,0,\mathsf z,\mathsf w)=\Lambda(\mathsf z\mathsf w^{-1},\mathsf B(a))$
and $K^{\mathrm c}(a,\pi/2,\mathsf z,\mathsf w)=\mathsf E_{\mathbb Z}(a(\mathsf z)a(\mathsf w)^{-1})$.
Then $K(a,\theta,\mathsf v)=K^{\mathrm c}(a,\theta,\mathsf z,1)\Lambda(\mathsf z,\mathsf Q)^{-1}$.

It is notable that loops which are already linear will remain constant but stabilized:
If $a(\mathsf z)=\Lambda(\mathsf z,\tilde Q)$ then
$K^{\mathrm c}(a,\theta,\mathsf z,\mathsf w)=\Diag(\ldots,\mathsf z\mathsf w^{-1}|
\Lambda(\mathsf z\mathsf w^{-1},\tilde Q),1\ldots)$, independently from $\theta$.
Similarly, rapidly decreasing perturbations of a linear loop will linearize
through rapidly decreasing perturbations of that linear loop.
\end{remark}
\begin{remark}
For a locally convex algebra $\mathfrak A$ we can define
\[K_0(\mathfrak A)=\pi_0^{\mathrm{smooth}}(\invol(\mathsf Q+\mathcal K_{\mathbb Z}(\mathfrak A))),\]
the smooth path components of the involutions, which are perturbations of $\mathsf Q$.
Similarly, one can  define
\[K_1(\mathfrak A)=\pi_0^{\mathrm{smooth}}(\unit(\mathcal K_{\mathbb Z}(\mathfrak A))).\]
Now $\mathsf B$, by this linearization argument, induces an isomorphism
\[\mathsf B_*:K_1(\mathfrak A[\mathsf z^{-1},\mathsf z]^{\mathrm{po}})\rightarrow K_0(\mathfrak A).\]
This is the ``hard part'' of Bott periodicity in the complex case, when geometric loops can be
represented by cyclic loops.
\end{remark}

\section{The contractibility of the pointed stable Toeplitz unit group} \label{sec:TP}

When we extend the stabilization procedure of Proposition \ref{prop:rapid} to Toeplitz algebras,
the symbol suddenly appears in the result:
\begin{proposition}\label{prop:toeploop} The continuous map
\[T:\mathcal T_{\mathbb N}(\mathfrak A)\times \mathrm S^1\rightarrow
\mathcal W_{\mathbb N}(\mathfrak A)+
\mathcal K_{\mathbb N}(\mathfrak A)[\mathsf v^{-1},\mathsf v]\subset
\mathcal T_{\mathbb N}(\mathfrak A)[\mathsf v^{-1},\mathsf v]\]
defined by
\[A,\theta\mapsto T(A,\theta,\mathsf v)=\delta_{t,1}a(\mathsf
v)\mathbf e_{0,0}+\delta_{t,-1}a(-\mathsf v)\mathbf e_{0,0}+
\mathsf V^{-1}\mathsf C(\theta)\mathsf VA\mathsf V^{-1}\mathsf
C(\theta)^\dag\mathsf V,\]
where $a=\sigma(A)$, is smooth in the variable $\theta$. It yields a family of homomorphisms of
$\mathcal T_{\mathbb N}(\mathfrak A)$ to
$\mathcal T_{\mathbb N}(\mathfrak A)[\mathsf v^{-1},\mathsf v]$.
The map leaves the symbol invariant.
For $\theta=0$,
\[T(A,0,\mathsf v)=A=\begin{bmatrix}
a_{11} & a_{12} &\cdots\\
a_{21} & a_{22} &\cdots\\
\vdots&\vdots&\ddots
\end{bmatrix};\]
but for $\theta=\pm\pi/2$,
\[T(A,\pm\pi/2,\mathsf v)=
\begin{bmatrix}
a(\pm\mathsf v)\\
&a_{11} & a_{12} &\cdots\\
&a_{21} & a_{22} &\cdots\\
&\vdots&\vdots&\ddots
\end{bmatrix}.\]
\begin{proof}
It follows from  direct inspection of the matrices in question.
\end{proof}
\end{proposition}
As a corollary we obtain
\begin{proposition}\label{prop:toepsymbol}
The continuous map
$Z:\unit(\mathcal T_{\mathbb N}(\mathfrak A))\times\mathrm S^1\rightarrow
\unit(\mathcal K_{\mathbb N}(\mathfrak A)[\mathsf v^{-1},\mathsf v]^{\pointed})$
defined by
$A,\theta\mapsto Z(A,\theta,\mathsf v)=T(A,\theta,\mathsf v)T(A,\theta,1)^{-1}$
is smooth in  $\theta$. For $\theta=0$  it yields
$Z(A,0,\mathsf v)=1_{\mathbb N},$
but for $\theta=\pm\pi/2$ it yields
$Z(A,\pm\pi/2,\mathsf v)=\mathsf E_{\mathbb N}(a(\pm \mathsf v)a(1)^{-1})$.

Consequently, the symbols $a(\mathsf z)$ of invertible Toeplitz algebra elements
are stably homotopic to constant loops $a(1)$.
If $\mathfrak A= \mathcal K_{\mathbb N}(\mathfrak S)$, then
(according to Corollary \ref{cor:infstable}) stable homotopy implies the existence of ordinary homotopies.
\qed
\end{proposition}

\begin{point}
Suppose that $Q$ is an involution, and $k\in\mathfrak A$.
We will use the shorthand notation
$k^+_Q=\frac12(k+QkQ)$,
$k^{++}_Q=\frac{1+Q}2k\frac{1+Q}2$,
$k^{+-}_Q=\frac{1+Q}2k\frac{1-Q}2$,
$k^{-+}_Q=\frac{1-Q}2k\frac{1+Q}2$.
Let us define
\[L(Q,k)= \mathsf W(\Lambda(\mathsf z,Q))k\mathsf W(\Lambda(\mathsf z^{-1},Q))=
\begin{bmatrix}
k^{++}_Q & k^{+-}_Q \\
k^{-+}_Q & k^{+}_Q & k^{+-}_Q \\
&k^{-+}_Q & k^{+}_Q &\ddots\\
&&\ddots&\ddots
\end{bmatrix}.\]
This is a homomorphism in $k$, and we can extend it to $\tilde k=1+k$ by
$L(Q,\tilde k)=1_{\mathbb N}+L(Q,k)$.
Notice that in this case, $\tilde kL(Q,\tilde k^{-1})$ has symbol
$\tilde k\Lambda(\mathsf z,Q)\tilde k^{-1}\Lambda(\mathsf z,Q)^{-1}$.
\end{point}
\begin{point}\label{po:constru}
Assume that $Q=\mathsf Q$ and $k\in \mathcal T_{\mathbb Z}(\mathfrak S)$.
Set
\[\widetilde{L}(k)=\left[\begin{array}{ccc|cccc}
\ddots&\ddots&&&&&\\
\ddots&k^{+}_Q &k^{-+}_Q&&&&\\
&k^{+-}_Q&k^{+}_Q&k^{-+}_Q\mathbf e_{00}&k^{-+}_Q\mathbf e_{10}&k^{-+}_Q\mathbf e_{20}&\cdots\\\hline
&&\mathbf e_{00}k^{+-}_Q&\mathbf e_{00}k^{++}_Q\mathbf e_{00}&
\mathbf e_{00}k^{++}_Q\mathbf e_{10}&\mathbf e_{00}k^{++}_Q\mathbf e_{20}&\cdots\\
&&\mathbf e_{01}k^{+-}_Q&\mathbf e_{01}k^{++}_Q\mathbf e_{00}&
\mathbf e_{01}k^{++}_Q\mathbf e_{10}&\mathbf e_{01}k^{++}_Q\mathbf e_{20}&\cdots\\
&&\mathbf e_{02}k^{+-}_Q&\mathbf e_{02}k^{++}_Q\mathbf e_{00}&
\mathbf e_{02}k^{++}_Q\mathbf e_{10}&\mathbf e_{02}k^{++}_Q\mathbf e_{20}&\cdots\\
&&\vdots&\vdots&\vdots&\vdots&\ddots
\end{array}\right].\]
What happens here, compared to $L(Q,k)$, is the following: We inflated the first row and column
to infinitely many rows and columns, and reordered the matrix.
Again, this is a homomorphism in $k$, and we can extend it to $\tilde k\in1_{\mathbb Z}+
\mathcal T_{\mathbb Z}(\mathfrak S)$ by taking $\widetilde{L}(\tilde k)=
1_{\mathbb Z\times\mathbb Z}+\widetilde{L}(k)$.
Assume now that $\tilde k\in\unit(\mathcal T_{\mathbb Z}(\mathfrak S))$, and
the symbol of its lower right quadrant is $a(\mathsf z)$.
Consider
\begin{multline}
\mathsf U(\Lambda(\mathsf z,\mathsf Q))\widetilde{L}(\tilde k)\Lambda(\tilde k^{-1},\mathsf Q_{\mathcal T_{\mathbb Z}(\mathfrak S)})
\mathsf U(\tilde k\Lambda(\mathsf z,\mathsf Q)^{-1}\tilde k^{-1})=
\notag\\
=\left[\begin{array}{ccc|cc}
\ddots&\ddots&&&\\
&\frac{1-\mathsf Q}2&\frac{1+\mathsf Q}2&&\\\hline
&&\frac{1-\mathsf Q}2&\frac{1+\mathsf Q}2&\\
&&&\ddots&\ddots
\end{array}\right]
\widetilde{L}(\tilde k)
\left[\begin{array}{ccc|ccc}
\ddots&\ddots&&\\
&\frac{1+\mathsf Q}2&\frac{1-\mathsf Q}2&&\\\hline
&&\tilde k\frac{1+\mathsf Q}2&\tilde k\frac{1-\mathsf Q}2&\\
&&&\ddots&\ddots
\end{array}\right]\tilde k^{-1}.
\notag
\end{multline}
From the observation
$\widetilde{L}(\tilde k)\Lambda(\tilde k^{-1},\mathsf Q_{\mathcal T_{\mathbb Z}(\mathfrak S)})
\in\unit(\mathcal T_{\mathbb Z}(\mathcal K_{\mathbb Z}(\mathfrak S)))$,
and a careful examination of the matrix product, we find that the resulting expression
is of shape
$\begin{bmatrix}
1_{\mathbb N}&\\&N(\tilde k)
\end{bmatrix}
\in\unit(\mathcal T_{\mathbb Z}(\mathcal K_{\mathbb Z}(\mathfrak S)))$;
where we introduced the notation $N(\tilde k)$ for the lower right
quadrant.
Then the component
 $N(\tilde k)\in\unit(\mathcal T_{\mathbb N}(\mathcal K_{\mathbb Z}(\mathfrak S)))$ has symbol
$\Lambda(\mathsf z,\mathsf Q)\mathsf E_{\mathbb Z}(a(\mathsf z))
\tilde k\Lambda(\mathsf z,\mathsf Q)^{-1}\tilde k^{-1}
=\mathsf E_{\mathbb Z}(a(\mathsf z))\Lambda(\mathsf z,\mathsf Q)
\tilde k\Lambda(\mathsf z,\mathsf Q)^{-1}\tilde k^{-1}$.
Let us set $G(a)=N(\mathsf U(a))$. This yields
\end{point}
\begin{proposition}\label{prop:symsection} The continuous map
$G:\unit(\mathfrak S[\mathsf z^{-1},\mathsf z])\rightarrow \unit(\mathcal T_{\mathbb N}
(\mathcal K_{\mathbb Z}(\mathfrak S)))$
is such that  the symbol of $G(a)$ is
$\mathsf E_{\mathbb Z}(a(\mathsf z))\Lambda(\mathsf z,\mathsf Q)
\mathsf U(a)\Lambda(\mathsf z,\mathsf Q)^{-1}\mathsf U(a)^{-1}.$
\qed
\end{proposition}

Now, according to Proposition \ref{prop:toepsymbol} and  Corollary \ref{cor:stablehom},
the mere existence of the map above implies that
the symbol $\mathsf E_{\mathbb Z}(a(\mathsf z))
\Lambda(\mathsf z,\mathsf Q)\mathsf U(a)\Lambda(\mathsf z,\mathsf Q)^{-1}\mathsf U(a)^{-1}$
is homotopic to $1_{\mathbb Z}$ for
$a(\mathsf z)\in\unit(\mathfrak A[\mathsf z^{-1},\mathsf z]^\pointed)$.
So, Proposition \ref{prop:symsection} can be considered as a reformulation of linearizability.

\begin{proposition}[$\Rightarrow$ Statement \ref{stat:toepcont}]
\label{prop:toepcont}
The unit group
$\unit(\,\mathcal T_{\mathbb N}(\mathcal K_{\mathbb Z}(\mathfrak S))^{\mathrm{po}})$
is smoothly contractible.
\begin{proof}
We prove the statement up to stabilization. Then stabilization can be removed according to
Corollary \ref{cor:stablehom}.

(a) First, consider any element
$A\in \unit(\,\mathcal T_{\mathbb N}(\mathcal K_{\mathbb Z}(\mathfrak S))^{\mathrm{po}})$.
According to Proposition \ref{prop:toepsymbol}, its symbol $a$ is (stably) homotopic to the
constant loop $1$.
Applying Proposition \ref{prop:symsection} to this homotopy, we see that it is sufficient to prove
that Toeplitz units with symbol
$\mathsf U(a)\Lambda(\mathsf z,\mathsf Q)\mathsf U(a)^{-1}\Lambda(\mathsf z,\mathsf Q)^{-1}$
can be contracted.

(b) Consider, again, $A$ as above. Let
\[Q=\begin{bmatrix}-1_{\mathbb Z}&\\&\mathsf Q \end{bmatrix}=
\left[\begin{array}{cc||cc}-1_{\mathbb N}&&&\\
&-1_{\mathbb N}&&\\\hline&&&\vspace{-11pt}\\\hline&&-1_{\mathbb N}&\\&&&1_{\mathbb N}\end{array}\right];\]
here the double lines show how we decompose this block matrix of $\mathbb Z\times\mathbb Z$ matrices
to a block matrix of $\mathbb N\times\mathbb N$ matrices.
Furthermore, let
\begin{multline}\notag
S(\theta)=\begin{bmatrix}s&&&-t\\&1&&\\&&1&\\t&&&s\end{bmatrix}
\begin{bmatrix}1_{\mathbb N}&&&&\\&1_{\mathbb N}&&\\&&\mathsf W(a^\top) & \mathsf Y(a^\top)\\
&&\mathsf Y(a) & \mathsf W(a)\end{bmatrix}\\
\begin{bmatrix}1_{\mathbb N}&&&&\\&\mathsf W(a^{-1})&
\mathsf Y(a^{-1})&\\&\mathsf Y((a^{-1})^\top)&\mathsf W((a^{-1})^\top) & \\
&& &1_{\mathbb N} \end{bmatrix}
\begin{bmatrix}1_{\mathbb N}&&&\\&A&&\\&&1_{\mathbb N}&\\&&&A^{-1}\end{bmatrix}
\begin{bmatrix}s&&&t\\&1&&\\&&1&\\-t&&&s\end{bmatrix}\\
\in\unit(\mathcal K_{\{1,2\}\times\mathbb Z}(\mathfrak S)),
 \end{multline}
and take $S(\theta)L(Q,S(\theta)^{-1})\in
\unit(\mathcal T_{\mathbb N}(\mathcal K_{\{1,2\}\times\mathbb Z}(\mathfrak S)))$.
This yields a homotopy between $S(0)L(Q,S(0)^{-1})$ and $S(\pi/2)L(Q,S(\pi/2)^{-1})$, which have symbols
\[S(0)\Lambda(\mathsf z,Q)S(0)^{-1}\Lambda(\mathsf z,Q)^{-1}=\begin{bmatrix}
1_{\mathbb Z}&\\&\mathsf U(a)\Lambda(\mathsf z,\mathsf Q)\mathsf U(a)^{-1}\Lambda(\mathsf z,\mathsf Q)^{-1}
\end{bmatrix}\] and
\[S(\pi/2)\Lambda(\mathsf z,Q)S(\pi/2)^{-1}\Lambda(\mathsf z,Q)^{-1}=\begin{bmatrix}
1_{\mathbb Z}&\\&1_{\mathbb Z}
\end{bmatrix},\]
respectively. Thus, Toeplitz units with symbol
$\Lambda(\mathsf z^{-1},\mathsf Q)\mathsf U(a)\Lambda(\mathsf z,\mathsf Q)\mathsf U(a)^{-1}$
can be deformed to Toeplitz units with trivial symbols.
According to part (a), it is sufficient to show that elements with
trivial symbol can be contracted.

(c) Now suppose that the symbol of a Toeplitz unit $A$ is  $1$. According to standard
stabilization arguments, we can assume that $A=\mathsf E_{\mathbb N}(\tilde k)$,
where $\tilde k=\begin{bmatrix}k_0&\\&1_{\mathbb N}\end{bmatrix}\in\unit(\mathcal K_{\mathbb Z}(\mathfrak S))$.
Let
$\tilde k(\theta)=
\begin{bmatrix}s&t\\-t&s\end{bmatrix}
\begin{bmatrix}\tilde k_0&\\&1_{\mathbb N}\end{bmatrix}
\begin{bmatrix}s&-t\\t&s\end{bmatrix}.$
Then $\tilde k(\theta)L(\mathsf Q,\tilde k(\theta))^{-1}$ yields a homotopy
between $\tilde k(0)L(\mathsf Q,\tilde k(0))^{-1}=\mathsf E_{\mathbb N}(\tilde k)=A$
and $\tilde k(\pi/2)L(\mathsf Q,\tilde k(\pi/2))^{-1}=1$.
\end{proof}
\end{proposition}
\begin{remark} (a)
If the  locally convex algebra $\mathfrak A$ is strong in the terminology of in \cite{L},
i. e. for all seminorm $p$ there is a seminorm $\tilde p$
such that $p(X_1\ldots X_n)\leq \tilde p(X_1)\ldots \tilde p( X_n)$ holds for all $n$, then
the proof can be much simplified: In that case, the associated algebras
are also strong, and the smooth homotopy lifting property holds for the symbol map.
Then, using Proposition \ref{prop:toepsymbol}, the proof of the
contractibility statement reduces to point (c) immediately, hence making
points (a) and (b), and the construction of \ref{po:constru} unnecessary.
One must note that Proposition \ref{prop:toepcont} above is much easier to prove than
Kuiper's Theorem about the contractibility of the unitary group.
See, e. g. \cite{WO}.

(b)
Stabilization was an important assumption in the previous statement. For example,
$\unit(\,\mathcal T_{\mathbb N}(\mathbb C)^{\mathrm{po}})$ is not contractible,
as it allows an extended, multiplicate determinant.
\end{remark}

\section{Possible modifications}\label{sec:otheralg}
Due to the nice properties of $U(a,\theta,\mathsf v)$,
Statement \ref{stat:finite} can be seen in a rather straightforward manner.
We remark that another such category is the category of Hilbert-Schmidt operators,
used by Pressley and Segal, \cite{PS}, Ch. 6.
Furthermore, with some extra work, the transformation parameter $\theta$
(i. e. $s$ and $t$ jointly)
can be replaced by $t$ entirely, extending the constructions as formal homotopies.

\section{Algebraically finite cyclic loops}\label{sec:practbott}
 A  practical disadvantage
of $\mathsf B(a)$ is that it is, in general, an infinite perturbation of $\mathsf Q$.
The exception is when $a\in\unit(\mathfrak A[\mathsf z^{-1},\mathsf z]^{\mathrm f})$, but this
is a rather restrictive condition from geometrical viewpoint.
We will show below that we can do well also in the case when $a$ can be represented by finite
loops but it is not in $\unit(\mathfrak A[\mathsf z^{-1},\mathsf z]^{\mathrm f})$.
\begin{point}
For $m\leq 0\leq n$,
we say that
the loop $a(\mathsf z)\in \unit(\mathfrak A[\mathsf z^{-1},\mathsf z])$ is an $L(m,n)$-finite loop if
$a(\mathsf z)=\sum_{m\leq j\leq n}a_j\mathsf z^k$.
A loop $a(\mathsf z)$ is an $R(m,n)$-finite loop if its inverse $a(\mathsf z)^{-1}$ is an
$L(-n,-m)$-finite loop.
For a finite sequence $F=\{(m_j,n_j)\}_{1\leq j\leq s}$, let
\[\mathfrak A_F=\{(a_s,\ldots,a_1)\,:\,a_j\in\unit(\mathfrak A[\mathsf z^{-1},\mathsf z])
 \text{ is $L(m_j,n_j)$ or $R(m_j,n_j)$-finite}\}.\]
We say that
$a\in\unit(\mathfrak A[\mathsf z^{-1},\mathsf z]$ is algebraically finite of type $F$ if
$a=a_s\ldots a_1$ for an element $(a_s,\ldots,a_1)\in\mathfrak A_F$.
\end{point}

\begin{point}
For $m\leq 0\leq n$,
we say that a matrix $A$ is an $L(m,n)$-perturbation of $A_0$ if
\[A=A_0+\sum_{m\leq i\leq n,\,j\in\mathbb Z}a_{i,j}\mathbf e_{i,j},\]
for $a_{i,j}$ chosen suitably. Similarly, we can define $R(m,n)$-perturbations by interchanging the
role of $i$ and $j$ in the expression above. An $(m,n)$-perturbation  is a matrix
which is both an  $L(m,n)$-perturbation and an $R(m,n)$-perturbation.

In what follows, we will always be concerned with perturbations of
$\Lambda(\mathsf s,\mathsf Q)$,
where $\mathsf s$ is equal to $1$, $-1$, or another formal variable $\mathsf v$.
Both $L(m,n)$-perturbations and $R(m,n)$-perturbations of
$\Lambda(\mathsf s,\mathsf Q)$ can be reduced to
$(m,n)$-perturbations by taking direct cut-offs of unwanted matrix elements:
\[\begin{bmatrix}\mathsf s&&\\L^-&M&L^+\\&&1 \end{bmatrix}\xrightarrow{\Red_{(m,n)}}
\begin{bmatrix}\mathsf s&&\\&M&\\&&1 \end{bmatrix}\xleftarrow{\Red_{(m,n)}}
\begin{bmatrix}\mathsf s&R^-&\\&M&\\&R^+&1 \end{bmatrix}.\]

The reduction $\Red_{(m,n)}$ is essentially taking away the
off-diagonal elements of a triangular block matrix (with respect to an appropriate ordering of the basis).
Sometimes it is practical to use the partial reduction $\Red_{(m,n)}^{[h]}=(1-h)\mathrm{Id} + h \Red_{(m,n)}$,
where $h$ is assumed to be a scalar variable.
Here the off-diagonal blocks are not taken away completely but multiplied by $1-h$.
It is useful to notice that (partial) reduction is a homomorphism as long as we restrict our
attention to matrices of appropriate block triangular shape. In particular,
invertible elements / involutions are reduced to invertible elements / involutions.
\end{point}
\begin{point}
The involutions $Q$ and $\bar Q$ are unipotently related  if
$\frac12(Q\bar Q+\bar QQ)=1$ holds.
In this case the expression $C(\bar Q,Q)=\frac{1+\bar QQ}2$ satisfies the identities
\[C(\bar Q,Q)^{-1}=C(Q,\bar Q)
\qquad\text{and}\qquad
C(\bar Q,Q)\,QC\,(\bar Q,Q)^{-1}=\bar Q.\]

More generally, $C(\bar Q,Q,h)=(1-h)1+h\frac{1+\bar QQ}2$ satisfies
the identities
\[C(\bar Q,Q,h)^{-1}=C(Q,\bar Q,h)
\quad\text{and}\quad
C(\bar Q,Q,h)\,QC\,(\bar Q,Q,h)^{-1}=(1-h)Q+h\bar Q.\]

This situation applies when, in the manner of the previous paragraph,
an involution $Q$ is reduced to an involution $\bar Q$.
\end{point}
\begin{lemma}\label{lem:bu}
If $a(\mathsf z)=\sum_{m\leq j\leq n}a_j\mathsf z^k$, $m\leq0\leq n$, then $U(a,\theta,\mathsf v)$ is an
$(m,n)$-perturbation of $\mathsf U(a)$.
\begin{proof}
This is immediate from \ref{po:udesc}.
\end{proof}
\end{lemma}
\begin{lemma}\label{lem:looppert}
Suppose that $A$ is an $(m',n')$-perturbation of $\Lambda(\mathsf s,\mathsf Q)$, where $m'\leq0\leq n'$.
Then we claim:

If $a$ is an $L(m,n)$- or $R(m,n)$-finite loop, then
$U(a,\theta_1,\mathsf v)A U(a,\theta_2,\mathsf w)^{-1}$ is an
$L(m+m',n+n')$- or $R(m+m',n+n')$-perturbation of $\Lambda(\mathsf s,\mathsf Q)$, respectively.
\begin{proof}
The $L$ case:
Let  $k>n+n'$ and $h=\mathsf s$ if
$k<m+m'$.
The special shape of the matrices implies
\[\mathbf e_k^\top U(a,\theta_1,\mathsf v)A=\biggl(\sum_{m\leq j\leq n}a_j\mathbf e_{k-j}^\top\biggr)A
=h\sum_{m\leq j\leq n} a_j\mathbf e_{k-j}^\top=h\mathbf e_k^\top U(a,\theta_2,\mathsf w),\]
from which $\mathbf e_k^\top U(a,\theta_1,\mathsf v)AU(a,\theta_2,\mathsf w)^{-1}=h\mathbf e_k^\top=
\mathbf e_k^\top\Lambda(\mathsf s,\mathsf Q)$.
This latter equality, which holds for appropriate $k$, is exactly the statement
of having an $L(m+m',n+n')$-perturbation of $\Lambda(\mathsf s,\mathsf Q)$.
The $R$ case is similar.
\end{proof}
\end{lemma}

\begin{point}
Next, we construct a linearization procedure which linearizes algebraically finite loops
into finite perturbations:
Let $F=\{(m_j,n_j)\}_{1\leq j\leq s}$ be a finiteness type,  $\tilde a=(a_s,\ldots, a_1)\in\mathfrak A_F$,
and $a=a_s\ldots a_1$.
Set $M_k=m_1+\ldots+m_k$, $N_k=n_1+\ldots+n_k$.
Let $|F|=(M_s,N_s)$.
Also, let  $\tilde a_k=(a_k,\ldots, a_1)$, with appropriate finiteness type $F_k$.
Then $|F_k|=(M_k,N_k)$.
We define
\[\mathsf B_F(\tilde a)=\Red_{|F_s|}\left(\mathsf U(a_s)\ldots
\Red_{|F_1|}\left(\mathsf U(a_1) \mathsf Q\mathsf U(a_1)^{-1}\right)
\ldots \mathsf U(a_s)^{-1}\right).\]
Then $\mathsf B_F(\tilde a)$ is an involution, and an $|F|$-perturbation of $\mathsf Q$.
More generally, let
\begin{multline}\notag
K_F^{\mathrm c}(\tilde a,\theta,\mathsf v,\mathsf w)=\Red_{|F_s|}\Bigl(U(a_s,\theta,\mathsf v)\ldots\\\ldots
\Red_{|F_1|}\Bigl(U(a_1,\theta,\mathsf v)\Lambda(\mathsf v\mathsf w^{-1},\mathsf Q)U(a_1,\theta,\mathsf w)^{-1}\Bigr)
\ldots U(a_s,\theta,\mathsf w)^{-1}\Bigr).
\end{multline}
Then, in particular,
$K_F^{\mathrm c}(\tilde a,0,\mathsf v,\mathsf w)=\Lambda(\mathsf v\mathsf w^{-1},\mathsf B_F(\tilde a))$,
and
$K_F^{\mathrm c}(\tilde a,\pi/2,\mathsf v,\mathsf w)=\mathsf E_{\mathbb Z}(a(\mathsf z))\cdot
\Lambda(\mathsf v\mathsf w^{-1},\mathsf Q)\mathsf E_{\mathbb Z}(a(\mathsf w)^{-1})$;
which  are immediate from the special shape of the matrices involved.
This  yields
\end{point}
\begin{proposition}
The continuous map
\[K_F:\mathfrak A_F\times\mathrm S^1\rightarrow
\unit(\mathcal K_{\mathbb Z}(\mathfrak A)[\mathsf v^{-1},\mathsf v]^{\pointed})\]
defined by
\[ \tilde a,\theta\mapsto
K_F(\tilde a,\theta,\mathsf v)=K_F^{\mathrm c}
(\tilde a,\theta,\mathsf v,1)\Lambda(\mathsf v,\mathsf Q)^{-1}
\]
is smooth in the variable $\theta$; and it is an $|F|$-perturbation of $1_{\mathbb Z}$.
Here
\[K_F(\tilde a,0,\mathsf v)=\Lambda(\mathsf v,\mathsf B_F(\tilde a))\Lambda(\mathsf v,\mathsf Q)^{-1},\qquad
K_F(\tilde a,\pi/2,\mathsf v)=
\mathsf E_{\mathbb Z}(a(\mathsf v)a(1)^{-1}).\]

In particular, as $\mathrm S^1$ is restricted to $[0,\pi/2]$, it yields a linearizing homotopy of $a(\mathsf z)a(1)^{-1}$
in the finite perturbation category. \qed
\end{proposition}
In the literature one finds comments about the possibly very large
size of the matrices used in linearizing homotopies. The result above, however,
shows the one can do reasonably well.

\begin{point} There is, however,  a closer analogy between the non-finite and the finite cases:
Let $Q_0=\mathsf Q$, and $Q_k=\Red_{|F_k|}( \mathsf U(a_k) Q_{k-1}\mathsf U(a_k)^{-1} )$ by recursion.
Then $Q_k=\mathsf B_{F_k}(\tilde a_k)$.
Using the notation $\prod_{i=1}^sx_i=x_n\ldots x_2x_1$,
let
\[{\mathsf U}_F(\tilde a)=\prod_{i=1}^s\frac{\mathsf U(a_i)+Q_i\mathsf U(a_i)Q_{i-1}}2
=\prod_{i=1}^sC\left(Q_i,\mathsf U(a_i)Q_{i-1}\mathsf U(a_i)^{-1}\right)\,\mathsf U(a_i).\]
According to our earlier observations,
\[\mathsf B_F(\tilde a)={\mathsf U}_F(\tilde a)\mathsf Q{\mathsf U}_F(\tilde a)^{-1}.\]

We also define
\[U_{F}(\tilde a,\theta, \mathsf v)=
\Red_{|F_s|}\left(U(a_s,\theta,\mathsf v)\ldots
\Red_{|F_1|}\left(U(a_1,\theta,\mathsf v)\mathsf U(a_1)^{-1}\right)
\ldots\mathsf U(a_s)^{-1}\right)
{\mathsf U}_F(\tilde a).\]
and
\[\widehat{\mathsf U}_F(\tilde a)=
\Red_{|F_s|}\left(\widehat{\mathsf U}(a_s)\ldots
\Red_{|F_1|}\left(\widehat{\mathsf U}(a_1)\mathsf U(a_1)^{-1}\right)
\ldots\mathsf U(a_s)^{-1}\right)
\mathsf U_F(\tilde a).\]

Then
$U_{F}(\tilde a,0, \mathsf v)=\mathsf U_F(\tilde a),$
which is trivial; and, analogously to the original situation,
$U_{F}(\tilde a,\pi/2, \mathsf v)=
\mathsf E_{\mathbb Z}(a(\mathsf z))\Lambda(\mathsf v,\mathsf Q)\widehat{\mathsf U}_F(\tilde a)
\Lambda(\mathsf v,\mathsf Q)^{-1},$
which follows from
$\Lambda(\mathsf v,\mathsf B_F(\tilde a))^{-1}=
{\mathsf U}_F(\tilde a) \Lambda(\mathsf v,\mathsf Q)^{-1}{\mathsf U}_F(\tilde a)^{-1}$
and the homomorphism property of reduction.
In fact,
\[K_F^{\mathrm c}(\tilde a,\theta,\mathsf v,\mathsf w)
=U_{F}(\tilde a,\theta, \mathsf v)
\Lambda(\mathsf v\mathsf w^{-1},\mathsf Q)U_{F}(\tilde a,\theta, \mathsf w)^{-1}\]
holds.
Again, this follows from
 $\Lambda(\mathsf v\mathsf w^{-1},\mathsf B_F(\tilde a))=
 {\mathsf U}_F(\tilde a)\Lambda(\mathsf v\mathsf w^{-1},\mathsf Q){\mathsf U}_F(\tilde a)^{-1}$
and the homomorphism property of reduction.
\end{point}

\begin{point}
The  constructions above can be expounded in order to show that
the linearizations $K$ and $K_F$ can nicely be deformed into each other:
Let
\[
\mathsf U_F(\tilde a,h)=
\prod_{k=1}^s (1-h) \mathsf U(a_k)+h\mathsf U_{F_k}(\tilde a_k)\mathsf U_{F_{k-1}}(\tilde a_{k-1})^{-1}\]
Here the product terms can also be written as
$ C(Q_k,\mathsf U(a_k)Q_{k-1}\mathsf U(a_k)^{-1},h)\mathsf U(a_k)$,
which makes invertibility clear.
Then
$\mathsf U_F(\tilde a,0)=\mathsf U(a)$, $\mathsf U_F(\tilde a,1)=\mathsf U_F(\tilde a)$.
Let
\[\mathsf B_F(\tilde a,h)={\mathsf U}_F(\tilde a,h)\mathsf Q {\mathsf U}_F(\tilde a,h)^{-1}.\]
Notice that
$\mathsf B_F(\tilde a,0)=\mathsf B(a)$, $\mathsf B_F(\tilde a,1)=\mathsf B_F(\tilde a).$
Let
\begin{multline}
U_{F}(\tilde a,h,\theta,\mathsf v)=\prod_{k=1}^s
\left((1-h)U(a_k,\theta,\mathsf v)+hU_{F_k}(\tilde a_k,\theta,\mathsf v)U_{F_{k-1}}(\tilde a_{k-1},\theta,\mathsf v)^{-1}\right)=\\
=\prod_{k=1}^s\Red_{|F_k|}^{[h]}\left(U(a_k,\theta,\mathsf v)\ldots
\Red_{|F_1|}\left(U(a_1,\theta,\mathsf v)\mathsf U(a_1)^{-1}\right)
\ldots\mathsf U(a_k)^{-1}\right)\notag\\\notag
C(Q_k,\mathsf U(a_k)Q_{k-1}\mathsf U(a_k)^{-1},h)\mathsf U(a_k)\\\notag
\Red_{|F_{k-1}|}\left(U(a_{k-1},\theta,\mathsf v)\ldots
\Red_{|F_1|}\left(U(a_1,\theta,\mathsf v)\mathsf U(a_1)^{-1}\right)
\ldots\mathsf U(a_{k-1})^{-1}\right)^{-1}.
\end{multline}
Again, the latter product form implies not only invertibility but that the
inverses of the product terms are  linear in $h$. In particular, it yields that  the inverse is
\[U_{F}(\tilde a,h,\theta,\mathsf v)^{-1}=\prod_{k=1}^s(1-h)U(a_k,\theta,\mathsf v)^{-1}+hU_{F_{k-1}}(\tilde a_{k-1},\theta,\mathsf v)
U_{F_k}(\tilde a_k,\theta,\mathsf v)^{-1}.\]
This also shows that $U_{F}(\tilde a,h,\theta,\mathsf v)^{-1}$  is polynomial in $h$.
We also define
\[\widehat{\mathsf U}_F(\tilde a,h)=\prod_{k=1}^s
\left((1-h)\widehat{\mathsf U}(a_k)+
h\widehat{\mathsf U}_{F_k}(\tilde a_k)\widehat{\mathsf U}_{F_{k-1}}(\tilde a_{k-1})^{-1} \right)\]
One can see that the identities
$U_{F}(\tilde a,h,0,\mathsf v)={\mathsf U}_F(\tilde a,h)$
and
\[U_{F}(\tilde a,h,\pi/2,\mathsf v)=
\mathsf E_{\mathbb Z}(a(\mathsf z))\Lambda(\mathsf v,\mathsf Q)
\widehat{\mathsf U}_F(\tilde a,h)
\Lambda(\mathsf v,\mathsf Q)^{-1}\]
hold. Furthermore,
$U_{F}(\tilde a,0,\theta,\mathsf v)=U(a,\theta,\mathsf v)$,
$U_{F}(\tilde a,1,\theta,\mathsf v)=U_{F}(\tilde a,\theta,\mathsf v)$,
and $\widehat{\mathsf U}_F(\tilde a,0)=\widehat{\mathsf U}(a)$,
$\widehat{\mathsf U}_F(\tilde a,1)=\widehat{\mathsf U}_F(\tilde a)$.
 We define
\[K_F^{\mathrm{ec}}(\tilde a,h,\theta,\mathsf v,\mathsf w)
=U_{F}(\tilde a,h,\theta, \mathsf v)
\Lambda(\mathsf v\mathsf w^{-1},\mathsf Q)U_{F}(\tilde a,h,\theta, \mathsf w)^{-1}.\]
From the earlier observations, the identities
$K_F^{\mathrm{ec}}(\tilde a,h,0,\mathsf v,\mathsf w)=\Lambda(\mathsf v\mathsf w^{-1}, \mathsf B_F(\tilde a,h))$
and
\[K_F^{\mathrm{ec}}(\tilde a,h,\pi/2,\mathsf v,\mathsf w)=
\mathsf E_{\mathbb Z}(a(\mathsf v))\Lambda(\mathsf v\mathsf w^{-1},\mathsf Q)\mathsf E_{\mathbb Z}(a(\mathsf w))^{-1}\]
follow. Furthermore,
$K_F^{\mathrm{ec}}(\tilde a,0,\theta,\mathsf v,\mathsf w)=
K^{\mathrm{c}}(a,\theta,\mathsf v,\mathsf w)$
and
$K_F^{\mathrm{ec}}(\tilde a,1,\theta,\mathsf v,\mathsf w)
=K_F^{\mathrm{c}}(\tilde a,\theta,\mathsf v,\mathsf w)$.

This yields
\end{point}
\begin{prop}[$\Rightarrow$ Statement \ref{stat:subfinite}]
The continuous map
\[K_F^{\mathrm{e}}:\mathfrak A_F\times \mathbb R\times\mathrm S^1\rightarrow \unit(\mathcal K_{\mathbb Z}
(\mathfrak A)[\mathsf v^{-1},\mathsf v]^{\mathrm{po}})\]
defined by
\[\tilde a,h,\theta\mapsto K_F^{\mathrm{e}}(\tilde a,h,\theta,\mathsf v)
=K_F^{\mathrm{ec}}(\tilde a,h,\theta,\mathsf v,1)\Lambda(\mathsf v,\mathsf Q)^{-1}
\]
is smooth in $\theta$ and polynomial in $h$. It has the properties

(i) $K_F^{\mathrm{e}}(\tilde a,0,\theta)=K(a,\theta)$;

(ii) $K_F^{\mathrm{e}}(\tilde a,1,\theta)=K_F(\tilde a,\theta)$;

(iii) $K_F^{\mathrm{e}}(\tilde a,h,0)=\Lambda(\mathsf v,\mathsf B_F(\tilde a,h)\Lambda(\mathsf v,\mathsf Q)^{-1}$;

(iv) $K_F^{\mathrm{e}}(\tilde a,h,\pi/2)=\mathsf E_{\mathbb Z}(a(\mathsf z)a(1)^{-1})$.

In particular, it connects the pullback homotopy $K|_{\mathfrak A_F}$ and homotopy $K_F$ through other linearizing homotopies.
\qed
\end{prop}

\end{document}